\documentclass{amsart}
\newtheorem{theorem}{Th\'eor\`eme}[section]
\newtheorem{lemma}[theorem]{Lemme}
\newtheorem{prop}[theorem]{Proposition}
\newtheorem{corollary}[theorem]{Corollaire}
\theoremstyle{definition}
\newtheorem{definition}[theorem]{D\'efinition}

\theoremstyle{remark}
\newtheorem{remark}[theorem]{Remarque}

\numberwithin{equation}{section}

\def\H{\mathbb H}
\def\C{\mathbb C}
\def\Q{\mathbb Q}
\def\P{\mathbb P}

\def\Z{\mathbb Z}

\def\NN{\mathcal N}
\def\HH{\mathcal H}
\def\LL{\mathcal L}
\def\PP{\mathcal P}

\def\n{\noindent}

\newcommand{\pt}{ \,^{p}\!\tau  }

\newcommand{\ilm}{  j_{!*}{\mathcal L}}
\newcommand{\ph}[2]{ \,^p\!{\mathcal H}^{#1}(#2)}

\title{\large Du  th\'eor\`eme
de d\'ecomposition \`a la Puret\'e locale }
\author{Fouad El Zein}
\address{Institut de Math\'ematiques de Jussieu,
Paris, France}
 \email{elzein@math.jussieu.fr}
 \author{ D\~ung Tr\'ang L\^e}
\address{Universit\'e d'Aix-Marseille 
LATP, UMR-CNRS 7353
39, rue Joliot-Curie
F-13453 Marseille Cedex 13, France}
 \email{ledt@ictp.it}

\date{F\'evrier  2013}
\begin{document}
\maketitle 
\begin{abstract} Une nouvelle d\'emonstration du th\'eor\`eme de d\'ecomposition est donn\'ee, en \'etablissant une relation  avec une version
du th\'eor\`eme de puret\'e locale de Deligne et Gabber  adapt\'ee aux vari\'et\'es
alg\'ebriques complexes.
\end{abstract}

\vskip.2in
\centerline {\large From   the Decomposition theorem to Local purity}

\begin{abstract}
A new proof of the decomposition theorem is established using a  relation with a version of  the local purity theorem of Deligne and Gabber  adapted  to complex algebraic varieties.
\end{abstract}

\section {Abridged version}
 Let $f: X \to V$ be a  projective morphism of complex algebraic
 varieties, $\tilde \LL$ a polarized
variation of Hodge structures  (VHS)  on an algebraic open dense
subset $j: (X-Y) \to X$ complement of  a normal crossing divisor (NCD) $Y$ in a non singular variety $X$, $\LL := \tilde \LL [m]$ is the complex
with sheaf cohomology $\tilde \LL$ in degree  $- m $ where $ m $ is the
dimension of $X$,
 and $j_{!*} \LL$ is
 the intermediate extension of $\LL$ (\cite{BBD}, prop. (2.1.11) p. 60, (2.2.7) p. 69, \cite{B}).
 
 We give a new proof of the decomposition theorem
 of the direct image of $\ilm$ into its perverse cohomology: 
$$
 Rf_*j_{!*}\LL \simeq \oplus_i {^p\!\HH}^i(Rf_*j_{!*}\LL)[-i]
 $$
   and the canonical decomposition of each  ${^p\!\HH}^i(Rf_*j_{!*}\LL)$  into a direct sum of intermediate  extensions (\cite{BBD}, theorems (4.3.1) p. 112, (5.4.10) p. 144, (6.2.5) p. 163, (6.2.10) p. 165)  and   (\cite{MI}), based  on  Deligne - Gabber's 
note on local purity (\cite{DG}, not published). The theorem is true for any varieties $X$ and $Y$, but we show that the proof can always be reduced to the case of a normal  crossing divisor $Y$ in a non-singular $X$. 

As the references show, this theorem has been established first for a scheme $X$ of finite type 
over an algebraically closed field $k$ of characteristic $p > 0$, for pure sheaves in the abelian sub-category of perverse sheaves of the category $D^b_c(X, \Q_l)$ for $l \not= p$, with respect to l-adic cohomology. The statements can be transposed
according to Deligne's dictionary  (\cite{IC}, sections 3 and 9), which is done in (\cite{BBD}, 6.2). 

We consider a stratification of $V$ and   start with the theorem over the open dense subset (the big strata) on which the restriction of $f$ is smooth, where   the decomposition applies in the smooth proper case by the original work of Deligne  \cite{SS}.

 In the next inductive  step we  aim to extend the decomposition to the union $S_{n-1}$ of the next strata of dimension $\leq n-1$ where $n$ is the dimension of $V$. We consider  a normal section $\NN$
 at a  point  $v \in S_{n-1}$:
 \begin{enumerate}
\item Assuming the decomposition on $ \NN - v$,  we prove
 an adapted version of  the local purity theorem \cite{DG}  to complex algebraic 
 varieties on $ \NN $ at $v$. 
\item Then we use this local purity to extend the decomposition along the strata. 
\end{enumerate}
We  prove  this last result in a second  note, then both results give a new proof of the decomposition theorem.  In fact, in the proof of local purity, we need  essentially to consider the case of points on the strata of dimension zero which  presents similarity with the isolated singularity case \cite{EI}. 
The  inductive step is for a simultaneous proof of both results: local purity and decomposition.

 The proof is also based  on an induction on the dimension of $X$.  We use, down on $V$, M. Artin's results 
of Lefschetz type with coefficients  in perverse sheaves on the complement of hyperplane sections \cite{BBD}, theorem (4.1.1). 

The use of Hodge theory is limited in this note to the case of  the complement $X- Z$ of a normal crossing divisor  $Z$,  in order to use  logarithmic complexes   (\cite{EII}  theorem (3.5) p. 30  and 3.3 p. 33).  

When 
$X - Z$ is the inverse image of an open subset of $V$, we  prove that \emph{ the perverse filtration 
of its cohomology with coefficients in $\ilm$ is compatible with the mixed Hodge structure (MHS)}. 
 
In the case  of a small ball $B_v$ with center $v$ with inverse image $B_{X_v} = f^{-1}(B_v)$, we suppose  that the fiber  $X_v := f^{-1}(v)$ is
a NCD by blowing up if necessary; similarly we may suppose that the union $X_v \cup Y$ is 
also a NCD in $X$ in order to have a MHS constructed by a logarithmic complex. We  assume 
the decomposition theorem true on $B_{X_v} - X_v$ by induction, from which we deduce the isomorphism:
$$  Gr^{\pt}_i \H^{r+i}(B_{X_v} - X_v, \ilm) \simeq \H^r(B_v- v, {^p\!\HH}^i(Rf_*j_{!*}\LL))$$
 so that we can derive an interpretation of the local purity as a condition on the weights of the MHS induced  on the  graded perverse cohomology  of $B_{X_v} - X_v $.

\section {Introduction}
Soit  $f: X \to V$ un  morphisme  projectif de 
 vari\'et\'es alg\'ebriques complexes, $\tilde \LL$ 
une variation de structures de Hodge  (VSH) polaris\'ee  sur un ouvert $ \Omega$ lisse de $X$, $j: \Omega \to X$, $\LL := \tilde \LL [m]$  vu comme un complexe de faisceaux \'egal \`a  $\tilde \LL$ en degr\'e  $ - m $ o\`u $ m  $
est la dimension de $X$ et 
nul en tout autre degr\'e. L'extension  interm\'ediaire de  $\LL$  
 (\cite{BBD}, prop. (2.1.11) p. 60, (2.2.7) p. 69, \cite{B})
 est not\'ee  $j_{!*} \LL$.
  
Le th\'eor\`eme de d\'ecomposition dans \cite{BBD}  
 \'etablit dans le cas g\'eom\'etrique, la d\'ecomposition dans la cat\'egorie d\'eriv\'ee 
$D^b_c (V, \Q)$: 
$$ Rf_*j_{!*}\LL \simeq \bigoplus_{i \in \Z} {^p\!\HH}^i(Rf_*j_{!*}\LL)[-i]$$
du complexe image directe d\'eriv\'ee en somme directe de ses cohomologies  perverses. D'autre part,  chaque terme se d\'ecompose naturellement cette fois en somme directe d'extensions interm\'ediaires   (\cite{BBD}, th\'eor\`emes (4.3.1) p. 112, (5.4.10) p. 144, 
 (6.2.5) p. 163, (6.2.10) p. 165). Pour  le cas d'une VSH admissible et polaris\'ee, voir \cite{MI}. 
 
 \smallskip
 Ce th\'eor\`eme a \'et\'e \'etabli d'abord pour les sch\'emas $X$ de  type fini
sur un corps  alg\'ebriquement clos  $k$ de caract\'eristique $p > 0$, pour des faisceaux purs dans la sous-cat\'egorie ab\'elienne des faisceaux pervers de la  cat\'egorie $D^b_c(X, \Q_l)$ d\'eriv\'ee des $\Q_l$-faiceaux constructibles pour $l \not= p$ \cite{BBD}.

 Les \'enonc\'es se transposent
selon le  dictionnaire \'etabli par  Deligne   (\cite{IC}, sections 3 et 9); ce qui est fait dans (\cite{BBD}, 6.2) o\`u l'on d\'ecouvre une technique qui permet de d\'eduire des r\'esultats comportant des \'enonc\'es g\'eom\'etriques en caract\'eristique $0$ \`a partir du cas correspondant en caract\'eristique $p > 0$.

\smallskip
 Soit $\pi: \tilde X \to X$ une d\'esingularisation adapt\'ee \`a $f$, la d\'emonstration pour $f$ se  d\'eduit facilement de celles pour 
 $\pi$ et $f\circ \pi$  ce qui nous ram\`ene \`a consid\'erer le  cas o\`u $X$ est 
 lisse. En fait, on peut se ramener au cas o\`u  $f: X \to V$ est une fibration en DCN sur les strates, au sens donn\'e plus bas, qui va nous permettre d'utiliser les complexes logarithmiques en th\'eorie de Hodge et le bon comportement de la perversit\'e en pr\'esence d'un diviseur de Cartier et qui explique la relative simplicit\'e de cette note.
 
   La preuve de la d\'ecomposition a \'et\'e pr\'ec\'ed\'ee d'une note  \cite{DG} \'etablissant 
   la notion de puret\'e locale en caract\'eristique $p > 0$. Nous allons \'etablir  
   dans cette note une 
 relation r\'ecurrente entre ces deux r\'esultats qui conduit \`a une d\'emonstration 
 simultan\'ee 
 des deux r\'esultats: puret\'e et  d\'ecomposition, \`a partir de:
  \begin{enumerate}
\item  La construction d'une structure de Hodge mixte (SHM) sur les ouverts 
compl\'ementaires d'un diviseur \`a croisements normaux (DCN) dans $X$, \'etablie \`a l'aide d'un complexe logarithmique  (\cite{EII} th\'eor\`eme (3.5) p. 30 et 3.3 p. 33).
\item   La notion de filtration 
perverse ${^p\!\tau}$ sur l'image directe $K := Rf_*j_{!*}\LL$ (\cite{BBD} prop. 1.3.3 p. 29)
construite dans la cat\'egorie d\'eriv\'ee 
$D^b_c (V, \Q)$ des faisceaux \`a cohomologie constructible.
\end{enumerate}
 Un r\'esultat essentiel de la note, consiste \`a d\'emontrer que  { \it la filtration perverse} induite sur la cohomologie d'une extension interm\'ediaire sur ouvert $\Omega $ qui est \`a la fois compl\'ementaire d'un DCN dans $X$ et  inverse par $f$ d'un ouvert sur $V$,  { \it est compatible avec la SHM}.
 
  La  { \it polarisation } des structures de Hodge (SH) en caract\'eristique $0$ dans les travaux   \cite{CKS, KI, K}, 
 va nous permettre  de 
{\it  simplifier le cas crucial de la d\'emonstration de la puret\'e locale}.

  Par cons\'equent,  on  supposera dans cette note  $X$ lisse, la VSH $\LL $ localement unipotente et polaris\'ee sur l'ouvert $\Omega := X - Y$  compl\'ementaire d'un diviseur \`a croisements normaux  $Y$ dans $X$.
  
Soit $Z$ un DCN dans $X$ tel que $Z \cup Y$ soit aussi un DCN  (en effectuant des \'eclatements si n\'ecessaire) et $j_Z: X-Z \to X$, alors  un 
 complexe  de  Hodge mixte (CHM) cohomologique sur $X$ sous-jacent  \`a 
$R(j_Z)_* j_Z^* j_{!*}\LL$ est d\'ecrit  \`a l'aide d'un  { \it complexe logarithmique } dans  \cite{EII, E} (ou d'un module de Hodge  diff\'erentiel dans  \cite{ MI});
 dans chaque cas  la construction utilise des r\'esultats de Kashiwara  \cite{K, KI} (ou l'appendice  \`a  \cite{MI}).
 
Pour tout entier $k \in \Z$, la  {\em  filtration perverse} $\pt$ sur la  cohomologie  globale ${\H}^k(X-Z, K)$ est d\'efinie par:
\begin{equation*}
    \pt_i {\mathbb H}^k(X-Z, K) \,  :=\, Im
    { \,\left\{ {\H}^k(X-Z,\pt_i K) \to {\H}^k(X-Z, K)
      \right\}}.
\end{equation*} 
 Il s'agit de v\'erifier qu'elle est compatible avec la SHM sur $X-Z$.
 
 \begin{prop}
 \label{mainp} i) Soit  $\LL $ une VSH polaris\'ee sur le compl\'ementaire d'un  diviseur \`a croisements normaux  $Y$ de $X$,
$ U$ un ouvert 
 alg\'ebrique de $ V$, tel que $f^{-1}(U)$ soit le
compl\'ement  d'un   diviseur \`a croisements normaux   $Z$  dans $X$ tel que   $Z \cup Y$ soit aussi un DCN, alors  la filtration  perverse sur la cohomologie
 $\H^r( f^{-1}(U),
j_{!*}\LL)$ est compatible avec la structure de Hodge mixte (SHM): les  sous-espaces ${^p\!\tau}_i$ sont des  sous-SHM de $\H^r( f^{-1}(U), j_{!*}\LL)$.

ii) Soient un point  $v \in V$  (resp. une boule   $B_v \subset V$ de centre $v$) d'image inverse $f^{-1}(v) := X_v$
(resp. un voisinage  $f^{-1}( B_v) := B_{X_v} $ de $X_v$).  On suppose   $ X_v $ et $ X_v \cup Y$
 des  DCN dans $X$, alors pour toute boule $B_v$ de rayon assez petit la 
cohomologie $ \H^r(B_{X_v} - X_v, j_{!*}\LL) $ est munie d'une SHM  compatible avec la filtration perverse.  
\end{prop}
 
 L'interpr\'etation de la puret\'e locale  \cite{WII, DG} se traduit dans ce cas par  l'\'enonc\'e principal suivant constituant
 le pas de l'induction  qui permet de d\'eduire ce r\'esultat en un point $v$ \`a partir de la d\'ecomposition en dehors de ce point. On peut supposer le point dans la strate de dimension z\'ero (sinon, si $v$ est sur une strate $S_l$ on peut appliquer  le r\'esultat sur une  section normale \`a $S_l$ en $v$ pour se r\'eduire \`a ce cas).

\begin{theorem}\label{main}   Avec les notations de \ref{mainp}, soit  $a$ le poids   de la  VSH polaris\'ee $\LL$ (apr\`es d\'ecalage) et supposons que  la fibre $X_v := f^{-1}(v)$ image inverse du point $v$
par le morphisme projectif  $f: X \to V$, soit un DCN tel que  $ X_v \cup Y$
soit aussi un  DCN dans $X$.\\
 Si  la restriction de $Rf_* j_{!*} \LL$ \`a une boule  \'epoint\'ee   $B_v - \{v\}$ assez petite  
 satisfait  le th\'eor\`eme de d\'ecomposition,
alors  la SHM de la cohomologie de $B_v - \{v\}$   induit sur les espaces ci-dessous des SHM 
 de poids $\omega$  
satisfaisant les in\'egalit\'es suivantes:
\begin{enumerate}
\item  $\omega >  a+r $ sur ${^p\tau}_{\leq r} \H^r(B_{X_v}-X_v, j_{!*}\LL) $,
 \item Dualement,  $\omega
\leq a+r $ sur $ \H^r(B_{X_v}-X_v,j_{!*}\LL) / {^p\!\tau}_{\leq r}
\H^r(B_{X_v} -X_v,j_{!*}\LL) $.
\end{enumerate}
  \end{theorem}
 Ce sont des conditions dites de puret\'e locale au point $v$ de $V$ et de semi-puret\'e sur $X$ en $X_v$.
 Notre but est d'\'etablir ces conditions  sur les  poids de la SHM induite sur le gradu\'e par rapport \`a $\pt$.
Les SHM sont toutes r\'ealis\'ees ici comme sous-quotient de SHM d\'efinie  en haut sur des ouverts de $X$ lisse.
Le th\'eor\`eme \ref{main} est le pas   d'une r\'ecurrence simultan\'ee avec la preuve du  th\'eor\`eme
de d\'ecomposition, r\'ecurrence sur la dimension des strates et la  dimension de $X$ \`a l'aide d'une r\'eduction \`a une section hyperplane g\'en\'erale sur $X$.

 La {\it r\'ecurrence}  qui sera expliqu\'ee dans une seconde note, d\'ebute par  le r\'esultat de Deligne sur l'ouvert $U$
de la grande strate lisse  de $V$ sur lequel  la restriction de $f$ est lisse et propre \cite{SS}. Un exemple simple est trait\'e  dans le  cas d'une singularit\'e isol\'ee \cite{EI}. 
 Les arguments importants de la preuve sont les suivants:
 
1)  L'hypoth\`ese  de   d\'ecomposition en dehors du point $v$, permet d'utiliser l'isomorphis\-me:  
\begin{equation*}
Gr^{^p\!\tau}_i
\H^r(B_{X_v}-X_v,j_{!*}\LL) \simeq
\H^{r-i}(B_v-v,{^p\!\HH}^i(Rf_*j_{!*}\LL)) = \H^{r}(B_v-v,Gr^{\pt}_i(Rf_*j_{!*}\LL)) 
\end{equation*}
  exprimant la d\'eg\'en\'erescense de la suite spectrale de Leray perverse en dehors de $v$, ce qui permet d'exploiter les r\'esultats d'Artin de type Lefschetz en bas sur $V$ \cite{BBD}.  
  
2) L'utilisation de la polarisation de la cohomologie d'intersection 
simplifie consid\'erablement le cas crucial que l'on rencontre en caract\'eristique  $p > 0$.

Dans une seconde note \cite{EL2} nous prouvons que la puret\'e locale entra\^ine la 
d\'ecomposition au point et par cons\'equent les deux r\'esultats 
d\'emontrent le th\'eor\`eme de d\'ecomposition.

\section{Preuve de la semi-puret\'e}

La SHM sur la cohomologie du compl\'ementaire d'un DCN   telle qu'elle est d\'evelopp\'ee 
dans \cite{EII} (th\'eor\`eme (3.4) et propri\'etes p. 31, section (3.3)  p. 33 et 
th\'eor\`eme 3.14) est sp\'ecialement adapt\'ee pour \'etudier la semi-puret\'e  sur  
les voisinages  tubulaires $B_v:=f^{-1}(B_v)$ de $X_v$ dans cette note.  Notamment l'isomorphisme de Thom-Gysin  \cite{EII} (3.3.7) sera souvent 
utilis\'e dans la d\'emonstration.
Par ailleurs, le morphisme $f: X \to V$ sera toujours une fibration par des DCN 
sur les strates au sens suivant:

\begin{definition}[fibration topologique par DCN sur les strates d'une stratification]
i)  Un morphisme $f: X \to V$ est une fibration topologique par DCN sur les strates d'une stratification ${\mathcal S}=(S_\alpha)$ de Thom-Whitney de V si  $X$ est lisse et les espaces
$V_l=\cup_{\dim S_\alpha \leq l} S_\alpha$ satisfont les propri\'t\'es suivantes:
\begin{enumerate}
\item  L'espace $X$ est lisse, et les  sous-espaces  $ X_{V_i} := f^{-1} ( V_i)$ sont vides ou des sous-DCN emboit\'es dans $X$.
\item (T) La restriction de $f$ \`a $X_S := f^{-1} (S)$ au-dessus de chaque  strate $ S$  de ${\mathcal S}$ est une fibration  topologique: 
 $f_{\vert}: X_{S} \to S$.
 \item  Pour tout point $v \in  V_i - V_{i-1}$, donc lisse dans $V_i$, soit $\NN_v$ une section normale en $v$ \`a $V_i$ en position g\'en\'erale, alors $f^{-1}(\NN_v)$ est lisse dans $X$ et intersecte le diviseur \`a croisements normaux  
   $X_{V_i} $ transversalement.
\end{enumerate}   
  Si ces deux derni\`eres assertions sont satisfaites,  alors 
  $X_{V_i} \cap f^{-1}(\NN_v) =  f^{-1} (v) $  est un DCN  dans  $f^{-1} (\NN_v)$.

\n  Nous dirons pour simplifier, que $X_{{V_i} - V_{i-1}}:= f^{-1} ( {V_i}- V_{i-1})$ est DCN relatif (\'eventuellement vide), et lorsque les donn\'ees sont claires, que le  morphisme $f$ (resp.  la stratification  $ \mathcal S $) est admissible. 

\n ii) La fibration  est adapt\'ee \`a un sous-espace  $Y$ dans $X$, ou \`a un syst\`eme  local  $\LL$ d\'efini sur le compl\'ementaire de 
 $Y$ dans $X$, si, de plus, $Y$ est un DCN et, pour tout $1\leq i\leq n$,
  la r\'eunion  des sous-espaces  
 $ X_{V_i} \cup Y$ sont des DCN relatifs  sur les strates  de $V$.
  \end{definition}
  
 Il existe  une stratification de Thom-Whitney sous-jacente aux donn\'ees de la fibration telle que 
les op\'erations cohomologiques classiques sur $\LL $ induisent des syst\`emes locaux sur les strates \cite{LT}.

\begin{remark}[r\'eduction au cas  d'une strate de dimension z\'ero]

\

 La notion de DCN relatif est essentielle, car elle permet de ramener l'\'etude du probl\`eme en un point $v$ d'une strate quelconque $S$ au cas d'un point $v$ d'une strate de dimension z\'ero dans la section transversale $\NN_v$ \`a $S$ en $v$.
  \end{remark}  
  
\begin{prop} 
Soient $f: X \to V $ un morphisme  projectif et  $Y$ un sous-espace alg\'ebrique ferm\'e strict contenant les singularit\'es de $X$, alors il existe un diagramme 
$X \xleftarrow{\pi'} X' \xrightarrow{f'}  V$ o\`u  $X'$ est une vari\'et\'e nonsinguli\`ere, $\pi'$ et  $f':= f \circ \pi'$ sont des fibrations par DCN  sur les strates, adapt\'ees \`a $Y':= \pi'^{-1}(Y)$. 

 De plus,  $\pi'$ est une modification  de $X$: il existe un ouvert  
 $\Omega \subset f(X) \subset V$ dense dans l'image de $f$, tel que $\pi'$ 
 induise un isomorphisme de 
 $f^{-1}(\Omega) - ( f^{-1}(\Omega) \cap Y) \simeq f'^{-1}(\Omega) - (f'^{-1}(\Omega) \cap Y')$, et que $f'^{-1}(\Omega) \cap Y'$ soit  un DCN relatif (dit horizontal, \'eventuellement vide).
  \end{prop}

Cette propri\'et\'e de $f$ est n\'ecessaire pour le raisonnement par  r\'ecurrence. 

\subsection{ La  filtration  perverse}

 La filtration perverse, bien que d\'efinie sur $V$, peut \^etre d\'ecrite directement sur $X$ d'apr\`es \cite{DM}. Cette description va nous permettre d'\'etablir la compatibilit\'e de la filtration perverse avec la SHM d'un ouvert alg\'ebrique \`a coefficients. Elle  s'adapte \`a la situation locale sur $V$ en un point $v$,  et semi-locale sur $X$ au voisinage de la fibre $X_v$ en $v$.
 \subsubsection{}
 Soit $U$ un ouvert  de $V$ et consid\'erons  deux  suites croissantes de  sous-vari\'et\'es ferm\'ees de $U$:
   $$H_*:  U = H_0 \supset H_{-1}  \supset \ldots  \supset H_{-n},\quad 
   W_*:  U = W_0 \supset W_{-1}  \supset \ldots  \supset W_{-n} $$
  On note  $h_i: (U- H_{-i} )\to U$  l'inclusion et  on  consid\`ere 
 un complexe $K \in D^b_c (U, \Q) $ \`a cohomologie constructible born\'ee.
  En reprenant les notations de \cite{DM} (remark 3.6.6), on d\'efinit  la filtration sur $\H^*(U,K)$: 
    \begin{equation}\label{delta}
\delta_p \H^*(U,K) := \hbox{Im}\{\oplus_{i-j = p} \H^*_{W_{-j}}(U,(h_i)_! h_i^* K) \to  \H^*(U,K) \} 
\end{equation}
   La filtration $\delta $  est l'aboutisse\-ment d'une suite spectrale.

 \begin{prop}[\cite{DM}, theorem 4.2.1]    
 Soit $U$ un ouvert quasi-projectif de $V$ et  $K \in D^b_c (U, \Q) $.  
Pour un choix convenable 
    des deux suites $H_*$ et $W_*$ mais assez g\'en\'eral  (relativement \`a un plongement affine dans un espace projectif et une stratification compatible), la  filtration $\delta$ (\ref{delta}) est \'egale \`a
    la filtration perverse $\pt$ \`a un d\'ecalage d'indices pr\`es.  
    \end{prop}
  
  La preuve,  bas\'ee sur  la notion 
      de r\'esolution par des faisceaux pervers  acycliques sauf en degr\'e $0$, utilise le r\'esultat suivant   \cite{Arap, DM} qui semble remonter \`a une lettre in\'edite de Deligne:
   \begin{lemma}
  Soient $U \subset \P^N$ un plongement affine, $H$ et $H'$ deux sections hyperplanes en position g\'en\'erale, les inclusions $j: (U-H\cap U) \to U$ et  $ j': (U-H'\cap U) \to U$.
        Pour tout faisceau pervers  $\PP$ sur $U$, on a:    $\H^r(U, j_!j^*j'_*(j')^! \PP) = 0 $ pour $r \neq 0$.
 \end{lemma}
    \begin{remark}
 i)   Le choix assez g\'en\'eral des
     deux sections hyperplanes  $H$ et $H'$  assure un isomorphisme $ j_!j^*j'_*(j')^! \PP) \simeq  j'_*(j')^! j_!j^* \PP) $.
     
       La preuve du lemme utilise le lemme d'Artin-Lefschetz faible appliqu\'e au faisceau pervers
     $j_!j^*j'_*(j')^! \PP$  sur $U$ et \`a son dual et le fait que  $U - H $ soit affine. 
     
 ii) Soient  un point  $v \in V$ et une boule $B_v \subset V $ de centre $v$, ce r\'esultat s'applique    lorsque l'on remplace $U$ par $B_v -v$ car le plongement est alors de Stein
 et un ouvert du type $ B_v - (H_{-i}\cap B_v) \subset B_v -v$ est aussi de Stein, en cons\'equence le lemme d'Artin-Lefschetz faible s'applique.
 \end{remark}
 On d\'eduit de la remarque:
  \begin{lemma}\label{Stein}
 Soient un point  $v \in V$ et une boule   $B_v \subset V$ de centre $v$. La  filtration $\delta$ (\ref{delta}) pour $U = B_v -v$ est \'egale \`a
    la filtration perverse $\pt$ \`a un d\'ecalage d'indices pr\`es.   
\end{lemma} 
  \begin{prop}
 Avec les notations de \ref{main}, pour tout ouvert quasi-projectif $U $ de $V$  tel que $f^{-1}(U)$ soit le compl\'ementaire d'un
  diviseur $D$ \`a croisements normaux dans $X$ et que $D \cup Y$  soit aussi un DCN, 
  la  filtration $\delta$,  et par cons\'equent la filtration perverse $\pt$  sur $\H^j( f^{-1}(U), \ilm)$, 
  en tout  degr\'e $j \in \Z$, sont compatibles avec la SHM.
  \end{prop}
  
 \noindent {\it Preuve de la proposition  \ref{mainp}}.
 Il s'agit de munir les termes de la filtration $\delta$ (\ref{delta}) de SHM.
\subsubsection{ Cas d'un ouvert affine}  Nous remarquons d'abord que le  choix  des suites   $H_*$ et $W_*$ dans $U$ \'etant assez g\'en\'eral, nous pouvons supposer que les suites sont induites par des suites dans $V$ not\'ees abusivement aussi $H_*$ et $W_*$  tel que:
  
   {\it Les diff\'erentes images inverses 
    $ H_{-i}' := f^{-1}(H_{-i})$ et  $ W_{-i}' := f^{-1} (W_{-j} ) $  soient non-singuli\`eres, que leurs 
     intersections le long de $ f^{-1}(H_{-i} \cap W_{-j} ) $ soient  transversales dans $X$, et qu'ils  coupent transversalement  
    un DCN fixe donn\'e dans $X$}.
    
     Ces conditions permettent d'utiliser un complexe logarithmique le long d'un DCN contenant toutes ces familles. 

{\it Isomorphismes de  Thom-Gysin.}  En g\'en\'eral, si  $H' $ est une sous-vari\'et\'e   non-singuli\`ere  de codimension  $r$ transversale 
  dans $X$ \`a  un diviseur \`a croisements normaux $ D \cup Y$ de sorte que l'on   puisse d\'efinir  un CHM  cohomologique sur $H'$ logarithmique en $ H' \cap (D \cup Y)$, on peut
   r\'ealiser   les {\it morphismes de  Thom-Gysin}  sous forme de  morphismes 
  de complexes logarithmiques  induisant des morphismes de  $SHM$ de type $(r,r)$  \cite{EII} (3.3.7):
  
  On a $ f^{-1}(U) = X - D$ et si on pose $j_D: f^{-1}(U)  \to X$, on trouve  pour $r=1$:
  : \begin{equation}
  i_{H'}^! (R(j_D)_* j_D^*j_{!*}\LL ) \simeq   i_{H' *} (R(j_{D\cap  H'})_* j_{D\cap  H'}^* i_{H'}^*j_{!*}\LL )[-2]
   \end{equation}
Le  cas  $r > 1$ est similaire au cas $r = 1$.
 On pose

\medskip
\n 
$K' :=  Rj_{D *}j_D^*j_{!*}\LL$,  \quad  $K = Rf_* K'$, \quad  $h_i: (V- H_{-i}) \to V$,  \quad 
$h'_i: (X- H'_{-i}) \to X$;

\medskip 
D'apr\`es  les triangles distingu\'es:
 $$ h'_{i !} (h'_i)^* K' \to K' \to   i_{H'_{-i} *}i_{H'_{-i}}^* K', \quad  i^!_{W'_{-j}} h'_{i !} (h'_i)^* K' \to  i^!_{W'_{-j}} K' 
 \stackrel{\varphi}{\rightarrow}   i^!_{W'_{-j}}  i_{H'_{-i} *}i_{H'_{-i}}^* K'$$ 
on peut  calculer la cohomologie:   $\H^*_{W_{-j}}(V,(h_i)_! h_i^* K)\simeq H^*_{W'_{-j}} (X, h'_{i !} (h'_i)^* K')$  \`a l'aide du  c\^one sur
$\varphi$ d\'ecal\'e, afin de mettre une SHM sur  ces groupes et donc sur les termes de la filtration $\delta$  (\ref{delta}).

 On va 
  transformer les 
termes du c\^one \`a l'aide d'isomorphismes de Thom-Gysin (d\'efinis
   par le drapeau  $W'_*$ dans $X$ chaque fois de codimension $1$ ou directement sur $W'_{-j} $ de codimension $j$). Ainsi par exemple le terme  $i^!_{W'_{-j}} i_{H'_{-i} *}i_{H'_{-i}}^* K'$ devient isomorphe au complexe  sur    $ W'_{-j} \cap H'_{-i}$  logarithmique en
   $D_{-j}\cap W'_{-j} \cap H'_{-i}$.
   
   On consid\`ere sur chaque $W'_{-k}$ le diviseur \`a croisements normaux  $D_{-k}:= D \cap W'_{-k}$, avec $j_k: (W'_{-k} - D_{-k} ) \to W'_{-k}$, $K'_k :=  Rj_{k *}j_k^*  i_{W'_{-k}}^* j_{!*} \LL$,    alors on a des isomorphismes de Thom-Gysin  
    $i^!_{W'_{-k-1}} K'_k[2] \simeq  K'_{k+1}$ en partant de $ K'_0 := K' $ pour $k=0$, jusqu'\`a
    $k = j-1$, de m\^eme on a une transformation similaire pour le terme $i^!_{W'_{-j}} i_{H'_{-i} *}i_{H'_{-i}}^* K'$ pour arriver au c\^one sur un morphisme $\phi_{i,j}: i_{W'_{-j}}^*K' [-2j] \to  i_{W'_{-j}}^*  i_{H'_{-i} *} i_{H'_{-i}}^*K' [-2j]$
     du complexe 
   $K'_j$   sur $ W'_{-j}$ logarithmique en $D_{-j} \cap W'_{-j}$  d\'ecal\'e, et de but  le complexe      $ i_{W'_{-k}}^* i_{H'_{-i} *}i_{H'_{-i}}^* j_{!*} \LL$   sur  $ W'_{-j} \cap H'_{-i}$  logarithmique en
   $D_{-j}\cap W'_{-j} \cap H'_{-i}$ d\'ecal\'e.\\

 \smallskip 
 {\it  En conclusion  la SHM sur $ \H^*_{W'_{-j}}(X,(h'_i)_! (h'_i)^* K') $ se  d\'eduit \`a 
   indices pr\`es, de  la SHM d\'efinie par le c\^one mixte du morphisme  de complexes 
   de Hodge mixte  logarithmiques  $\phi_{i,j}$}. 

\subsubsection{ Cas du voisinage local \' epoint\'e $B_v^* := B_v - v$ d'un point $v$ de $V$}
  La preuve  est similaire \`a celle qui  pr\'ec\`ede.
Pour calculer la cohomologie du voisinage tubulaire moins la fibre centrale
  $B_{X_v}^* = B_{X_v} - X_v$   \`a coefficients  dans $\ilm$, on consid\`ere  le morphisme
 compos\'e suivant 
\begin{equation} I:  R i_{X_v}^{!}\ilm
\rightarrow j_{!*}\ilm \rightarrow   i_{X_v}^* \ilm
 \end{equation}
que l'on appelle morphisme d'intersection $I$. Soit  $j_{X_v}: (X-X_v) \to X$,  alors $I$ s'inscrit 
dans un triangle $R i_{X_v}^{!}\ilm
 \rightarrow   i_{X_v}^* \ilm \rightarrow  i_{X_v}^* Rj_{X_v*} j_{X_v}^* \ilm  \xrightarrow {[1]}$.
Par  cons\'equent  la cohomologie de  $B_{X_v}^*$  \`a coefficients  dans $\ilm$ 
est celle du c\^one sur $I$.
\begin{definition}
   La fibre $X_v$ \'etant suppos\'ee un DCN,  les termes   $R i_{X_v}^{!}\ilm
$ et $   i_{X_v}^* \ilm$  sont  munis 
   d'une structure de complexe de Hodge mixte cohomologique,
 et par cons\'equent   la cohomologie de  $B_{X_v}^*$ se trouve d'une munie de la SHM
  d\'efinie par le  c\^one  mixte sur $I$.
   \end{definition}
Pour montrer que cette SHM est  compatible avec la filtration perverse, on 
cherche \`a appliquer la 
formule (\ref{delta})  qui  caract\'erise  la filtration perverse d'apr\`es \ref{Stein}.

En reprenant en abusant les notations  pr\'ec\'edentes avec $B_{X_v}$ au lieu de $X$, $X_v$ au lieu de $D$, $j_{!*}\LL_{\vert B_{X_v}}$ au lieu de $j_{!*}\LL$, deux familles $H'_{-i}$ et $W'_{-j}$ dans $B_{X_v}$,  $h'_i: (B_{X_v} - H'_{-i}) \to B_{X_v}$  et $ i^!_{W'_{-j}} j_{!*}\LL_{\vert B^*_{X_v}} 
 \stackrel{\varphi}{\rightarrow}   i^!_{W'_{-j}}  i_{H'_{-i} *}i_{H'_{-i}}^* j_{!*}\LL_{\vert B^*_{X_v}}$,
 le groupe de cohomologie    $ \H^*_{W'_j}(B_{X_v}^*,(h'_i)_! (h'_i)^* \ilm)$ se d\'eduit  alors d'un  c\^one   sur $\varphi$  et la SHM,  dont on la munit, 
d'un double  c\^one mixte une fois sur $I$ et une fois  sur  $\varphi$ \`a partir du carr\'e:
\[
\begin{array}{ccc}
 i_{X_v \cap W'_{-j}}^! i_{W'_{-j}}^*\ilm  &  \stackrel{I_j}{\rightarrow} &  i_{X_v \cap W'_{-j}}^* \ilm \\
  \downarrow \varphi&   &\downarrow \varphi  \\
 i_{X_v \cap W'_{-j}}^!  i_{ W'_{-j}}^*  i_{H'_{-i} *} i_{H'_{-i}}^*\ilm & \stackrel{I^i_j}{\rightarrow}  &    i_{X_v \cap W'_{-j}}^*  i_{H'_{-i} *} i_{H'_{-i}}^*\ilm
\end{array}
\]
En conclusion, la filtration $\delta$ sur $\H^*(B_{X_v}^*, \ilm)$ coincide \`a indices pr\`es 
avec la filtration perverse et de plus elle est r\'ealis\'ee par des sous-SHM.

\subsection{Preuve du th\'eor\`eme \ref{main} }
C'est le r\'esultat principal  de cette note. Par hypoth\`ese on suppose   la d\'ecomposition 
au-dessus de $V - v$ satisfaite, ce qui se traduit par l'isomorphisme:
\begin{equation}\label{1}
Gr^{\pt}_i \H^r(B_{X_v}-  X_v, \ilm)\simeq
\H^{r-i}(B_v-\{v\},\ph{i}{ Rf_*\ilm }).
 \end{equation}
La preuve  utilise  l'isomorphisme pr\'ec\'edent  afin d'exploiter 
  le  th\'eor\`eme d'annulation  d'Artin du type section  hyperplane de  Lefschetz 
  \`a coefficients  dans $^pD_c^{\leq 0}(V)$  (\cite{BBD}  4.1),
  sur un ouvert de Stein (resp. affine) en bas dans $V$   (le r\'esultat \'etant local, on  peut supposer $V$ projective). 
  
En supposant $X_v$ un DCN, le terme  de  gauche de (\ref{1}) se trouve  muni d'une SHM induite.  Pour la commodit\'e de l'argumentation  on raisonne sur le terme de droite, puis on interpr\`ete le r\'esultat
\`a gauche, autrement dit, par abus de notation  on transporte la SHM sur le terme de droite (on ne sait pas encore   que cette structure est canonique mais c'est en fait notre but ultime qui sera l'objet d'une autre \'etude  \`a comparer avec la structure d\'efinie dans la th\'eorie des modules diff\'erentiels de Hodge).

Dans la suite pour simplifier les notations nous posons $K= Rf_*\ilm $, alors il faut  prouver  que le poids $ \omega$ de cette   
$SHM$   satisfait les  in\'egalit\'es:
\begin{equation*}
\begin{split}
\omega > a+i+j  \,\, \hbox{sur}\,\, \H^j(B_v-v,\ph{i}{K}) \,\, \hbox{si}\,\,  j \geq 0,\\
 \hbox{et}  \, \, \omega \leq a+i+j  \, \,  \hbox{sur}\, \,  \H^j(B_v-v,\ph{i}{K}) \, \,  \hbox{si} \, \,   j \, \, 
\leq -1
 \end{split}
  \end{equation*}
  Les deux cas cit\'es sont duaux.
La preuve pour $j > 0 $ s'obtient par  une {\it r\'ecurrence simple  sur la dimension appliqu\'ee \`a une section hyperplane g\'en\'erale de $V$, l'applica\-tion du  th\'eor\`eme d'Artin-Lefschetz et 
le morphisme de Gysin pour la section hyperplane g\'en\'erale}. 
C'est une adaptation  de celle de  \cite{DG} au cas transcendant et diff\`ere donc par l'utilisation de la polarisation et repose sur le fait que si  $\ilm $ 
est pure de poids $a$, les poids des  SHM   sur 
$\H^*(X-Z,\ilm)$ et  $\H_Z^*(X,\ilm)$  sont  $ \geq a$ (resp. $\leq a$ sur  les SHM duales 
$\H^*(Z,i_Z^* \ilm)$,  et $ \H_c^*(X-Z, \ilm)$).

\medskip
 
 \noindent 1) {\it Preuve pour $j >0$, similaire \`a \cite{DG}.} Soit $H$ une section  hyperplane  g\'en\'erale de 
 $V$ contenant le  point $v$ et soit $H_v = B_v \cap H$. 
En particulier $H$ est normalement plong\'ee en dehors du point $v$ et son
transform\'e  strict   $H'$ (ne contenant pas $X_v$) dans  $X$ 
est  transverse aux sous-espaces   
 $X_l:= f^{-1} (S_l)$ inverse des strates $S_l$ de $V$; en particulier les  troncations perverses
 commutent avec  la restriction \`a $H-\{v\}$. 
 
 On consid\`ere  la suite exacte de groupes de  cohomologie \`a coefficients dans $\ph{i}{ Rf_*\ilm }$: 
 \begin{equation*}
 {\H}^{j-2} (H_v- \{v\})(-1) \simeq {\H}^j_{H_v- \{v\}} (B_v-
\{v\})  \rightarrow  {\H}^j
(B_v- \{v\}) \rightarrow \H^j (B_v - H_v )
\end{equation*}
   Si $B_v $ est de dimension $1$, $B_v -v$  est de Stein, sinon $B_v -H_v$ est de Stein et sa cohomologie  $\H^j (B_v - H_v,\ph{i}{K})$  s'annule pour $j \geq 1$
d'apr\`es  le 
th\'eor\`eme d'annulation sur un ouvert de Stein
  appliqu\'e \`a $ \ph{i}{K}$  dans $ {^pD}_c^{\leq 0}V$. 
 On en d\'eduit que le morphisme de Thom-Gysin: 
 $${\H}^{j-2} (H_v- \{v\}, \ph{i}{K})(-1) \xrightarrow{G_j} {\H}^j
(B_v- \{v\}, \ph{i}{K})   $$
 est surjectif pour $j = 1$, et  un isomorphisme
pour $j > 1$. 
  
   Soit $f'$  la restriction de $f$ 
  \`a $H'$ et $K':= Rf'_*( j_{!*}\LL_{\vert H'})$. Par transversalit\'e 
en dehors de $v$, on a:
$\H^{j-2} (H_v- \{v\}, \ph{i}{K})  \simeq \H^{j-2} (H_v- \{v\}, \ph{i+1}{K'[-1]}$.
 En raisonnant par {\it r\'ecurrence } sur la dimension de $V$,  l'hypoth\`ese de d\'ecomposition sur $H' - (X_v \cap H') $ au-dessus de $H-\{v\}$ s'applique \`a  la restriction  $ j_{!*}(\LL_{\vert H'}[-1])$ de poids $a-1$  sur $H'$. 

Donc,   le terme de la SHM sur  $\H^{j-2} (H_v- \{v\}, \ph{i}{K})$ s'annule pour le poids 
  $\omega' > a-1+ i+1+j-2 = a+i+j-2$ 
  et  d'apr\`es la surjectivit\'e du morphisme de Thom-Gysin
  $G_j$,  on d\'eduit que le terme de la SHM sur $B_v- \{v\}$ s'annule pour $\omega = \omega' +2 > a-1+ i+1+j-2 = a+i+j$ (en ajoutant $2$, \`a cause  du twist des SHM).
  
\medskip
\noindent 2) {\it Preuve pour $j = 0$ et $\omega  \geq  a+i+j$.}   
Soient   $H_1$ une  section  hyperplane g\'en\'erale ne 
contenant pas  $v$,  $k_v: (V-v) \to V$ et $i_v:\{v\}\to V$ les immersions canoniques. 
On consid\`ere  la suite exacte:
\begin{equation*}
 \H^0(V - H_1,R k_{v *} k_v^*  \ph{i}{K})
   \rightarrow H^0( i_v^* R k_{v *} k_v^* \ph{i}{K})
      \rightarrow  \H^1( V-H_1, R k_{v !} k_v^* \ph{i}{K})
\end{equation*}
o\`u on utilise l'\'ecriture suivante: 
 $\H_c^1( V-H_1,  \ph{i}{K}) = 
\H^1( V-H_1, R k_{v !}  k_v^* \ph{i}{K})$  avec   $R k_{v !}  k_v^* \ph{i}{K} \in {^pD}_c^{\leq 0}V$ pour d\'eduire que la cohomologie $\H^1$  de l'espace affine $V-H_1$ 
  s'annule  d'apr\`es Artin.

Le terme \`a gauche $ {\H}^0(V - H_1,R k_{v *} k_v^*  \ph{i}{K})$ est de poids $\omega \geq a+i$ car il s'agit d'un ouvert.
Donc, le terme du milieu $ \H^0(B_v-\{v\}, \ph{i}{K}) \simeq H^0( i_v^* R k_{v *}
k_v^* \ph{i}{K})$, isomorphe \`a $Gr^{\pt}_i \H^i(X-X_v- , \ilm)$, est de poids $\omega \geq a+i$.

\medskip
\noindent {\it  Le cas crucial.}  Le th\'eor\`eme  (\ref{main}) est \'etabli pour $j > 0$ et par dualit\'e  pour $j < -1$, pour $j=0$ si $\omega  \geq  a+i+j$ et par dualit\'e  pour $j = -1$ 
si  $\omega  \leq  a+i$.  Il reste \`a   l'\'etablir pour  $\omega  <  a+i$ et $j = -1$, 
c-\`a-d il faut d\'emontrer l'annulation des termes:
\begin{equation}
Gr^W_{a+r } Gr^{\pt}_r
\H^{r-1}(B_{X_v}-X_v, \ilm )  \simeq Gr^W_{a+r} \H^{-1}(B_v-v,
\ph{r}{K}) 
\end{equation}
 pour tout $r \in \Z$. La preuve se fait en plusieurs \'etapes. Elle utilise l'id\'ee que la cohomologie 
 de $B_v-v$ s'inscrit dans deux suites exactes issues des triangles: 
$$
 i_v^* R k_{v !} k_v^* \ilm \to i_v^* \ilm \to i_v^* R k_{v *} k_v^* \ilm, \quad    R k_v^!  \ilm \to i_v^* \ilm \to i_v^* R k_{v *} k_v^* \ilm 
$$
 Soient $X^* = X -  X_v$ et $B_{X_v}^* = B_{X_v}  -  X_v$, alors  on en d\'eduit le diagramme  \`a coefficients dans $\ilm$
$$\begin{array}{ccccccccccc}
& \H^{r }_{X_v}(X)&&\to&&\H^{r }(X_v )&&&&\H^{r +1}_{X_v}(X)&\\
 &&\searrow&& \nearrow &   & \searrow&&\nearrow &&\\
(3.5) \qquad &  && \H^{r}(X)&&& &
  \H^{r }(B_{X_v}^*)&&\qquad &\\
 & &\nearrow&& \searrow &   & \nearrow&&\searrow &\\
  & \H_c^{r }(X^*)&&\to&& \H^{r }(X^*) &&&&\H_c^{r+1 }(X^*)&
   \end{array}$$
ce qui  permet de  r\'ealiser la cohomologie  de $B_{X_v}^*$  
\`a l'aide  de l'un  des deux c\^ones  mixtes
$$ R\Gamma _{X_v} ( X, \ilm)  \to R\Gamma ( X_v, \ilm), \quad  R\Gamma _c( X^*, \ilm)  \to R\Gamma ( X^*, \ilm).$$  On peut comparer les deux \`a l'aide du 
 c\^one  mixte form\'e par les termes \`a gauche du diagramme:
$ R\Gamma _c( X^*, \ilm)  \oplus  R\Gamma _{X_v} ( X, \ilm)  \to R\Gamma ( X, \ilm)$.

  On introduit d'abord deux morphismes  $\partial_X$ et $\alpha_X$
 et on en d\'eduit des morphismes injectifs sur des espaces gradu\'es convenables.
 
 \begin{lemma}\label{i1} 
   Le morphisme de  connexion suivant est injectif: 
$$Gr^W_{a+r}Gr^{\pt}_r\H^{r -1}(B_{X_v}^*, \ilm )
\stackrel{\partial_X}{\rightarrow} Gr^W_{a+r}Gr^{\pt}_r\H_c^{r
}(X^*, \ilm ).$$
 \end{lemma}
 \noindent On consid\'ere
la suite exacte longue:
$$
 \H^{-1}(V -v, \ph{r}{K}) \to \H^{-1}( B_v -\{v\}, \ph{r}{K})
\stackrel{\partial_V}{\rightarrow} H_c^0(V-\{v\},\ph{r}{K}) \to
 $$
Il suffit  de montrer  que:  $ \H^{-1}( V-\{v\}, \ph{r}{K})$   est pur de poids
 $ a+r-1$ ou par  dualit\'e $ H^1_c( V-\{v\},
\ph{r}{K})$ est pur de poids $ a+r+1$.  On consid\`ere la section hyperplane $H_1$ de $V$  ci-dessus  ne passant pas par $v$ et 
 la suite exacte:

\smallskip
$\H_{H_1}^1 (V,{R k}_{v
!}k_v^* \ph{r}{K}) \xrightarrow{\varphi} \H^1 (V,{R jk}_{v !}k_v^*
\ph{r}{K})
\rightarrow \H^1 (V - H_1 , R k_{v !}k_v^* \ph{r}{K}) $

\smallskip
\n  o\`u  $\varphi$ est surjectif car le dernier terme  s'annule; 
or l'espace 

\smallskip
  $\H^{-1} (H_1,{R k}_{v !}k_v^* \ph{r}{K})(-1)\simeq \H_{H_1}^1 (V,{R k}_{v
!}k_v^* \ph{r}{K}) $ 

\smallskip 
\n a une SH pure de poids $a+r+1$  
(la SH est d\'efinie en tant que sous-quotient de cohomologie   sur $H'_1 $ lisse dans $X$, 
elle est pure par r\'ecurrence sur  $H_1$).
 
 \smallskip
 Maintenant, on consid\`ere la suite exacte longue:
 $$ \H^{r-1}(X_v,\ilm) \rightarrow \H_c^r(X^*,  \ilm)
\stackrel{\alpha_X}{\rightarrow} \H^{r}(X, \ilm)$$

 \begin{lemma} \label{i2} Le morphisme suivant induit par $\alpha_X$  est  injectif:\\
 $ Gr^W_{a+r}\H_c^{r }(X^*, \ilm )
\stackrel {\alpha_X} {\rightarrow}
 Gr^W_{a+r}\H^{r}(X, \ilm )\simeq \H^r(X, \ilm ) $.
 \end{lemma}
Le poids de   $\H^{r-1}(X_v,\ilm)  $
est $< a+ r $ car $X_v$ est ferm\'e, donc le morphisme $Gr^W_{a+r}
\alpha_X$ est injectif. On peut n\'egliger de prendre \`a droite le gradu\'e pour $W$, car $ \H^{r}(X,
 \ilm)$ est pur de poids $a+r$.
  \subsubsection{}
Pour d\'emontrer $Gr^W_{a+r}Gr^{\pt}_r\H^{r -1}(B_{X_v}^*, \ilm ) = 0 $, une difficult\'e provient du fait que $\partial_X$ intervient modulo  $Gr^{\pt}_r$  alors que l'on est sans hypoth\`ese sur  
  la filtration $\pt$ en $v$ pour utiliser  $\alpha_X$ sur le gradu\'e de  
  $\pt$. En fait, on  consid\`ere plut\^ot le morphisme compos\'e   $ \gamma := \alpha_X  \circ \partial_X$ suivant (avec un abus de notation pour $ \partial_X$ modulo  $\pt$ ou non):
$$   Gr^W_{a+r}\H^{r -1}(B_{X_v}^*, \ilm )
\stackrel {\partial_X} {\rightarrow}
 Gr^W_{a+r}\H_c^{r }(X^*, \ilm ) \stackrel {\alpha_X} {\rightarrow}
 Gr^W_{a+r}\H^{r}(X, \ilm ), \quad  \gamma := \alpha_X  \circ \partial_X $$ 
   et on  d\'emontre que toute   classe $ {\overline u } \in Gr^W_{a+r}Gr^{\pt}_r\H^{r -1}(B_{X_v}^*, \ilm ) $ 
 peut  se relever en un \'el\'ement  $ u \in
{\pt}_{\leq r} Gr^W_{a+r} \H^{r -1}(B_{X_v}^*, \ilm )$, de telle fa\c{c}on que son image $ \gamma ( u) $ soit nulle dans $ Gr^W_{a+r}\H^{r}(X, \ilm )\simeq \H^{r}(X, \ilm)$. De la relation  $\alpha_X  ( \partial_X (u)) = 0$, on  d\'eduit par le lemmme (\ref{i2}) que $\partial_X (u) = 0$,  puis
 $\partial_X ({\overline u }) =  \hbox {cl}( {\partial_X (u)}) = 0$ et
 finalement par le lemmme (\ref{i1}) on a: ${\overline u }= 0$.
 
\begin{lemma}[Lemme principal]  Soit $ \gamma := \alpha_X  \circ \partial_X$.
Toute classe:
$$ {\overline u  } \in Gr^{\pt}_{ r} Gr^W_{a+r} \H^{r
-1}(B_{X_v}^*, \ilm )$$ 
 est  repr\'esent\'ee par un \'el\'ement $ u \in{\pt}_{\leq r} Gr^W_{a+r} \H^{r -1}(B_{X_v}^*, \ilm )$ 
d'image $ \gamma ( u) $ nulle dans $\H^{r}(X, \ilm)$
\end{lemma}
\begin{corollary} $Gr^{\pt}_{ r} Gr^W_{a+r}
\H^{r -1}(B_{X_v}^*, \ilm ) = 0$.
\end{corollary}
\noindent {\it L'id\'ee est de choisir $u$  tel  que l'\'el\'ement $\gamma (u) $ 
soit primitif}, afin 
 d'utiliser la polarisation sur   $\H^{r}(X, \ilm)$ pour d\'eduire  que $ \gamma (u)$ est nul.\\
 La preuve  se subdivise en plusieurs  \'etapes.\\
    Soit $H$ une section hyperplane g\'en\'erale de  $
X $ transverse \`a toutes les  strates de  $Y$. On pose $H_v:=  X_v
\cap H$, $B_{H_v}  := H \cap B_{X_v} $, et on consid\`ere  les immersions
 $j': (H - H \cap Y) \to H$, $k':(f (H) - \{v\}) \to f(H)$, la restriction  $ \ilm_{\vert H}$  de
 $ {\ilm}$ \`a $H$ muni d'un isomorphisme canonique
  $\rho:  \ilm_{\vert H} \xrightarrow{\sim}  \ilm_{\vert H}$  car $H$ est  transversale aux strates.
   
   \begin{lemma}  Soient $j_H: (X-H) \to X$ et  $(\LL_{(X-H)})_! := (j_H)_!j_H^*\ilm$,
alors on peut repr\'esenter ${\overline u }$  par un \'el\'ement $ u \in
{\pt}_{ r} Gr^W_{a+r} \H^{r -1}(B_{X_v}-X_v, \ilm ) $ \'egal \`a l'image d'un 
\'el\'ement
 $u! \in {\pt}_{ r} Gr^W_{a+r} \H^{r
-1}(B_{X_v}-X_v, (\LL_{(X-H)})_!) $.
\end{lemma}
\noindent On a: \\
 $\H^{r -1}(B_{H_v}-H_v, \ilm ) \simeq
    \H^{r}(B_{H_v}-H_v, j'_{!*} {\LL}_{\vert
   H}[-1])$  o\`u  $ {\LL}_{\vert  H}[-1]$  est de poids  $a -1$, d'o\`u par r\'ecurrence sur $H$:
   $$G^W_{a+r}Gr^{\pt}_{ r}\H^{r -1}(B_{H_v}-H_v, \ilm )
    \simeq G^W_{a+r}Gr^{\pt}_{ r+1} \H^{r}(B_{H_v}-H_v,j'_{!*} {\LL}_{\vert
   H}[-1]) = 0. $$
   De la suite  exacte courte  $(\LL_{(X-H)})_! \to \ilm \to (i_H)_* i_H^* \ilm  $,
on obtient la suite  exacte longue:
\begin{equation}\label{2}
\H^{r -1}(B_{X_v}-X_v, (\LL_{X-H})_!) \to   \H^{r -1}(B_{X_v}-X_v,
\ilm )\to   \H^{r -1}(B_{H_v}-H_v, \ilm )
\end{equation}
La suite  gradu\'ee associ\'ee pour $W$ et $\pt$, c'est \`a dire la suite obtenue  de  (\ref{2}) 
en appliquant $G^W_{a+r}Gr^{\pt}_{ r} $, 
  reste exacte  par hypoth\`ese de d\'ecomposition en dehors de 
  $X_v$ et son terme de droite  $G^W_{a+r}Gr^{\pt}_{ r}\H^{r-1}(B_{H_v}-H_v, \ilm )$ 
  s'annule par l'hypoth\`ese  de r\'ecurrence sur la semi-purit\'e sur
 $H$. Par cons\'equent  l'\'el\'ement ${{\overline u } }$  dans le terme du milieu 
 est l'image d'un \'el\'ement ${\overline u }!$ \`a gauche.
  On choisit un  repr\'esentant $u!\in G^W_{a+r}{\pt}_{ r}
\H^{r -1}(B_{X_v}-X_v, (\LL_{(X-H)})_! )$ de ${\overline u }!$, et on prend pour  $u$ l'image de
$u! $ dans  $G^W_{a+r}{\pt}_{ r} \H^{r -1}(B_{X_v}-X_v, \ilm )$.\\

\noindent {\it Preuve du lemme principal}. \\
\n 1) L'\'el\'ement {\it $ \gamma (u)$ est primitif}.  En effet l'\'el\'ement  
  $ u \in G^W_{a+r}{\pt}_{ r} \H^{r -1}(B_{X_v}-X_v, \ilm )$ d\'efini comme image de 
$u!\in G^W_{a+r}{\pt}_{ r} \H^{r -1}(B_{X_v}-X_v,  (\LL_{(X-H)})_!  )$ 
est l'\'el\'ement annonc\'e dans le lemme,  puisque sa classe est
${\overline u }$ et de plus  sa restriction  $\rho (\gamma (u)) \in G^W_{a+r}\H^{r}(H, \ilm) = \H^{r}(H, \ilm)$ 
\`a  $H$ s'annule, et par cons\'equent $ \gamma (u) \in \H^{r}(X, \ilm)$
est un \'el\'ement primitif. \\
2)  Il reste \`a prouver que $ \gamma (u) $ s'annule. On utilise  le diagramme  commutatif:
$$\begin{array}{ccc} \H^{r -1}(B_{X_v}-X_v, \ilm )&\stackrel{\partial_X}{\rightarrow}&
 \H_c^{r}(X - X_v, \ilm)  \\
 \downarrow \partial &   & \downarrow \alpha_X \\
   \H_{X_v}^{r}(X, \ilm) &\stackrel{A}{\rightarrow}&
   \H^{r}(X, \ilm) \end{array}$$
   qui r\'ealise $ \gamma (u)= \alpha_X  \circ \partial_X(u) = A \circ \partial
   (u)$,
    comme l'image d'un \'el\'ement $ \partial (u) \in Gr^W_{a+r}\H_{X_v}^{r}(X,
   \ilm)$, et on consid\`ere le diagramme
$$\begin{array}{ccccc} &&
Gr^W_{a+r} \H^{r}(X, \ilm) && \\
 &  A {\nearrow} & &{A^*}{\searrow}& \\
  Gr^W_{a+r} \H_{X_v}^{r}(X, \ilm) && \stackrel{ I}{\longrightarrow}& &
  Gr^W_{a+r} \H^{r}(X_v, \ilm) \end{array}$$
et son dual
$$\begin{array}{ccccc} &&
Gr^W_{a-r} \H^{-r}(X, \ilm) && \\
 &   {\nearrow} A & &{A^*}{\searrow}& \\
  Gr^W_{a-r} \H_{X_v}^{-r}(X, \ilm) && \stackrel{ I}{\longrightarrow}& &
  Gr^W_{a-r} \H^{-r}(X_v, \ilm) \end{array}$$
 Soit
$\eta $ le cup-produit  avec la classe d'une section hyperplane 
de $X$,  alors la polarisation $Q$ sur la partie  primitive de $ \H^{r}(X, \ilm)$
est d\'efinie, \`a l'aide du  produit de  dualit\'e de  Poincar\'e $P $ et  de l'op\'erateur $C $ de 
Weil,  par la formule:  $Q(v,b) := P (Cv , \eta^r ( {\overline b}))$.
 On a aussi le produit non-d\'eg\'en\'er\'e  d\'efini par
dualit\'e:
$${ P_v:  Gr^W_{a+r} \H_{X_v}^{r}(X,
\ilm)\otimes Gr^W_{a-r} \H^{-r}(X_v, \ilm) \to \C.}$$
La dualit\'e entre   $ A$ et $ A^*$ est d\'efinie pour tout
$b \in Gr^W_{a+r} \H_{X_v}^{r}(X, \ilm)$  et tout  $ c \in
Gr^W_{a-r}\H^{- r}(X, \ilm)$ par la formule :
$\; P ( A b ,  c ) = P_v( b, A^* c ).$

\smallskip
\noindent  A l'\'el\'ement $ u \in {\pt}_{\leq- r} Gr^W_{a-r} \H^{-r
-1}(B_{X_v}-X_v, \ilm )$ correspond  une image  par le
morphisme de connexion  $ \partial u := u' \in Gr^W_{a-r}
\H_{X_v}^{-r}(X, \ilm)$ tel que $\gamma (u) = A (u' )$.
Soit  $C $ l'op\'erateur de Weil d\'efini par la $SH$ sur $
Gr^W_{a-r}\H_{X_v}^{-r}(X, \ilm)$, alors:
$$P (C.
A u', \eta^r A \overline {u' }) = P_v(C u', A^*
\circ A (\eta^r \overline {u'})) = P(C.u', \eta^r I (\overline {u'})).$$
Or:  $I (u')= A^*\circ A \circ  \partial (u) = A^*\circ  \alpha_X  \circ \partial_X (u) = 0 $ 
car $A^*\circ  \alpha_X = 0$; on  en d\'eduit $\, P (C.
A u', \eta^r A \overline {u' } ) = 0$, donc $Au' = 0$ par
polarisation, ce qui termine la preuve.

\bibliographystyle{amsalpha}

\bibliographystyle{amsalpha}

\end{document}